\documentclass[a4paper,10pt]{article}
\usepackage{amssymb,amsmath,amsthm, amsbsy }

\newtheorem{theorem}{Theorem}[section]

\newtheorem{lemma}[theorem]{Lemma}
\newtheorem{proposition}[theorem]{Proposition}

\theoremstyle{definition}

\newtheorem{remark}[theorem]{Remark}

\begin{document}

\title{\textbf{\Large Twisted waveguide with a Neumann window}}

 \author{Philippe Briet$$,
Hiba Hammedi$$
 \footnote {briet@univ-tln.fr, hammedi@univ-tln.fr.}}
\date{\small  
\emph{
\begin{quote}
Aix-Marseille Universit\' e, CNRS, CPT, UMR 7332,
Case 907, 13288 Marseille  and 
Universit\' e de Toulon, CNRS, CPT, UMR 7332
83957, La Garde, France;\\
\end{quote}
}
January 2016}

\maketitle

\centerline{{\it Dedicated to Pavel Exner on the occasion of his 70th birthday} }

\begin{abstract}
This paper  is concerned with the  study of  the
existence/non-existence  of the  discrete spectrum  of the Laplace
operator  on  a domain  of $\mathbb R ^3$ which consists in a
twisted tube. This operator is defined  by means  of  mixed boundary
conditions.  Here we impose Neumann Boundary conditions on  a
bounded open  subset of the boundary  of the domain  (the Neumann
window) and Dirichlet boundary conditions  elsewhere.
\end{abstract}

{ \bf classification}
Primary 81Q10; Secondary 47F05.

\bigskip

{ \bf keywords}
Waveguide, mixed boundary conditions, twisting.

\maketitle

\section{Introduction}
In this work, we would like to study the influence of a geometric
twisting   on
   trapped modes which occur in  certain waveguides.
 Here the waveguide consists  in a  straight tubular   domain $\Omega_{0} :=  \mathbb{R}\times\omega$   having a  Neumann   window  on its boundary $\partial \Omega_0$.  \\
 The cross section $\omega$  is supposed to be an open bounded connected subset of $\mathbb{R}^2$
 of diameter $d>0$ which is not rotationally invariant.  Moreover  $\omega$  is supposed  to  have  smooth boundary $\partial\omega$.

 It can be shown that the  Laplace operator associated to such a straight tube has bound states \cite{H}.\\

Let us introduce  some notations.
 Denote  by ${\mathcal N}$ the Neumann window.  It   is  an  open bounded subset of the boundary  $\partial \Omega_0$. Let  ${\mathcal D}$ be  its complement set in   $\partial \Omega_0$.
  When  ${\mathcal N}$ is an annulus    of size $l>0$  we will denote it  by,
  $${\mathcal A}_a(l):=I_a(l) \times \partial \omega, I_a(l):=(a,l+a), a \in \mathbb R.$$
Consider first  the   self-adjoint operator $H_0^{\mathcal N}$
associated  to  the  following quadratic form. Let $D(Q^{\mathcal N})
=\{\psi\in\mathcal{H}^1(\Omega_0)\ \mid  \psi_{\lceil{\mathcal D}}=0
 \}$ and  for $\psi \in D(Q^{\mathcal N}) $,
$$ Q^{\mathcal N}(\psi)= \int_{\Omega_{0}}  \vert \nabla \psi \vert^2 dx$$
i.e. the Laplace operator  defined on $\Omega_0$   with Neumann
boundary conditions (NBC) on ${\mathcal N}$ and Dirichlet
boundary conditions  (DBC)  on ${\mathcal D}$ \cite{DK,K}. \\

It is  actually shown in the Section 2  of this paper  that if
${\mathcal N} $  contains an annulus of   size $l$ large enough then
$H^{\mathcal N}_0$ has at least one discrete eigenvalue. In fact it
is  proved in   \cite{H} that this holds true if  ${\mathcal N} $
contains an annulus of   any  size $l>0$.

The question we are interested in is  the following: is it possible
that the  discrete spectrum  of $H^{\mathcal N}_0$  disappears  when
we apply a geometric twisting on the guide? This question is
motivated by  the results of \cite{EKK, Kre} where it is shown
that this phenomenon  occurs  in   some bent tubes when they are
subjected to a   twisting defined from
  an angle function  $\theta$ having   a derivative $\dot \theta$ with a compact support.    In this paper  we consider
   the  situation     described above  which is very   different from  the one of \cite{EKK, Kre}.\\
    Let us   now define   the twisting \cite{BKRS,ExK}. Choose  $\theta  \in C_c^1(\mathbb R)$   and introduce  the  diffeomorphism
\begin{eqnarray}\label{L}
  \mathcal{L}:\Omega_0 &\longrightarrow & \mathbb{R}^3 \\
   \nonumber (s,t_{2},t_{3})&\longmapsto& \Big(s,t_2\cos\theta(s)-t_3\sin\theta(s),t_2\sin\theta(s)+t_3\cos\theta(s)\Big).
\end{eqnarray}
The twisted tube is  given  by
$\Omega_\theta:=\mathcal{L}(\Omega_0)$.
  Let $D(Q^{\mathcal N}_\theta ) =\{\psi\in\mathcal{H}^1(\Omega_\theta)\ \mid \ \psi_{ \lceil{ {\mathcal L}({\mathcal D})}}=0\}  $ and   consider the following quadratic form
\begin{equation}\label{q1}
Q^{\mathcal N}_\theta(\psi) :=  \int_{\Omega_{\theta}}  \vert \nabla
\psi \vert^2 dx, \; \psi\in
 D(Q^{\mathcal N}_\theta).
\end{equation}
Through  unitary equivalence, we then have to consider
\begin{eqnarray}\label{qatheta}
q^{\mathcal N}_\theta(\psi)&:=& Q^{\mathcal N}_\theta(\psi o
\mathcal{L}^{-1})=\parallel\nabla
'\psi\parallel^2+\parallel\partial_s\psi+\dot{\theta}\partial_\tau\psi\parallel^2,\quad
\end{eqnarray}
 $ \psi \in  D(q_\theta^{\mathcal N}):=\{\psi\in\mathcal{H}^1(\Omega_0)\ \mid\psi_{ \lceil{ {\mathcal D}}}=0 \}$
and where
\begin{equation}\label{not}
\nabla ':={}^t \left(\partial_{t_{2}},\partial_{t_{3}}\right), \quad
\partial_{\tau}:=t_{2} \partial_{t_3}-t_{3} \partial_{t_2}.
\end{equation}
Denote by $H_\theta^{\mathcal N}$ the associated  self-adjoint
operator. It is defined as follows (see \cite{DK,K}). Let $
  D(H_\theta^{\mathcal N}) = \{\psi\in
D(q_\theta^{\mathcal N}), \quad H^{\mathcal N}_\theta \psi\in
L^2(\Omega_0)\quad \frac{\partial\psi}{\partial n}{ \lceil_{ {\mathcal N}}}=0  \} $
 with
\begin{equation}  \label{Htheta}
H_\theta^{\mathcal N} \psi =
(-\Delta_\omega-(\dot{\theta}\partial_\tau+\partial_s)^2) \psi, \;
\end{equation}
where the transverse Laplacian $\Delta_{\omega}:=
\partial_{t_2}^{2}+
\partial_{t_3}^{2}$. If $ {\mathcal N}={\mathcal A}_a (l), l>0 $, we  will denote  these forms  respectively as    $Q^l_\theta, q^l_\theta$ and the  corresponding operator as $H^l_\theta $ and
if ${\mathcal N}= \emptyset$ we denote the  associated operator by $
H_\theta$.

Then  the main result of this paper is
\begin{theorem}\label{resultat}
i) Under  conditions stated above on $\omega$ and $\theta$, there
exists $l_{min}:=l_{min}(\omega,d)>0$ such as   if for some $ a \in
\mathbb R $ and $l>l_{min}$,  ${\mathcal N}  \supset {\mathcal
A}_a(l)$ then
\begin{equation}
\sigma_d(H_{\theta}^{\mathcal N} )\neq\emptyset.
\end{equation}

ii) Suppose  $\theta $  is  a non zero function    satisfying the same conditions as in i) and     has   a  bounded second
derivative. Then there exists $d_{max}:=d_{max}(\theta,\omega)>0$
such that for all $0<d \leq d_{max}$ there exists
$l_{max}:=l_{max}(\omega, d, \theta)$ such as for all $0<l \leq l_{max}$,
if ${\mathcal N} \subset {\mathcal A}_a(l) $ and supp$(\dot \theta)
\cap I_a(l) =\emptyset $ for some $ a \in \mathbb R $ then
\begin{equation}
\sigma_{d}(H_\theta^{\mathcal N} )=\emptyset.
\end{equation}
\end{theorem}

 Roughly speaking  this result implies that   for $d$ small enough,  the discrete spectrum      disappears  when  the   width  of the Neumann window    decreases.

Let us describe briefly the content of the paper. In the Section 2
we give the proof  of the Theorem \ref{resultat} i). The section 3 is
devoted to the proof of the second part of the Theorem
\ref{resultat}, this proof needs several steps. In particular  we
first establish a local Hardy inequality. This    allows  us to
reduce the problem to the analysis of a one dimensional Schr\"odinger
operator   from which the Theorem  \ref{resultat} ii) follows.
Finally in the  Appendix of the paper we give partial results we use
in previous sections.

\section{Existence of bound states }

 First we  prove the following. Denote by $ E_1, E_2, ....$ the eigenvalues  (transverse modes)
of the  Laplacian $-\Delta_\omega$ defined on $\mathrm{L}^2(\omega)$  with
DBC on $ \partial\omega$.  Let  $\chi_1, \chi_2,...$  be  the
associated eigenfunctions. Then we have
\begin{proposition}\label{EssSpec}
$\sigma_{ess}(H^{\mathcal N}_\theta)=[E_1,\infty)$.
\end{proposition}

\proof  We know that   $\sigma_{}(H_\theta)=[E_1,\infty)$ see e.g.
\cite{Briet-2009}. But  by usual arguments \cite{ RS}, $H^{\mathcal
N}_\theta \leq   H_\theta $, then
\begin{equation}\label{firstinclusion}
[{E}_1,\infty) \subset \sigma_{ess}(H^{\mathcal N}_\theta).
 \end{equation}
Let  $a' \in \mathbb R$  and $l'>0$ large enough  such that
${\mathcal N} \subset {\mathcal A}_{a'}(l')=I_{a'}(l') \times \partial \omega$ and
$ supp(\dot \theta) \subset I_{a'}(l')$. Let  $\tilde H^{l'}_\theta$
be the operator defined as in \eqref{Htheta} but with additional
Neumann boundary conditions  on $ \{a'\}\times \omega \cup
\{a'+l'\}\times \omega$. So  $H^{\mathcal N}_\theta \geq   \tilde
H^{l'}_\theta $ and then $  \sigma_{ess}(H^{\mathcal N}_\theta)
\subset \sigma_{ess}(\tilde H^{l'}_\theta)$ \cite{ RS}.

But  $\tilde H^{l'}_\theta= \tilde H_i \oplus \tilde H_e$. The  interior operator $\tilde H_i $ is the corresponding operator defined on $ L^2( I_{a'}(l') \times \omega) $ with NBC on $ \{a'\}\times \omega \cup  \{a'+l'\}\times \omega \cup {\cal N}$ and DBC elsewhere on $ {\mathcal A}_{a'}(l')$. By general arguments of \cite{RS} it has only discrete spectrum consequently $\sigma_{ess}(\tilde H^{l'}_\theta)= \sigma_{ess}(\tilde H_e)$.\\
Now  the exterior operator $ \tilde H_e $ is defined on $ L^2
((-\infty,a') \times \omega \cup (a'+l', \infty) \times \omega) $
with DBC on $(-\infty,a') \times \partial \omega  \cup (a'+l',
\infty)  \times \partial \omega$ and NBC on  $ \{a'\}\times \omega
\cup  \{a'+l'\}\times \omega $. Since $\theta =0$ for $x<a'$ and $x>
a'+l'$, it is easy to see that
$$\tilde H_e =  \underset{n \geq1 }{\oplus } (  -\partial^2 + E_n) (\chi_n, .) \chi_n.$$
Hence $\sigma(\tilde H_e) =\sigma_{ess}(\tilde H_e)=[E_1, +
\infty)$. \qed


\medskip

The  Theorem \ref{resultat} $i)$ follows from

\begin{proposition}\label{thmexistence} Under  conditions of the Theorem \ref{resultat} i),
there exists $l_{min}:=l_{min}(\omega,d)>0$ such as for all
$l>l_{min}$ we have
\begin{equation}
\sigma_d(H_{\theta}^l)\neq\emptyset.
\end{equation}
\end{proposition}

\proof
  Let $\varphi_{l,a}$ be
 the  following function
 \begin{center}
 $\varphi_{l,a}(s):=\left\{
              \begin{array}{ll}
                \frac{10}{l}(s-a), & \hbox{ on $ [a,a+\frac{l}{10})$;} \\
                1, & \hbox{ on $[a+\frac{l}{10},a+\frac{9l}{10})$;} \\
                -\frac{10}{l}(s-l-a), & \hbox{ on $[a+\frac{9l}{10},a+l)$;} \\
                0, & \hbox{elsewhere.}
              \end{array}
            \right.$
\end{center}
 It is  easy to see that   $\varphi_{l,a} \in D(q_\theta ^l)$ and $\parallel\varphi_{l,a}\parallel^2=\frac{13l}{15}\mid\omega\mid$.
Let us calculate
\begin{equation}\label{eq353}
q^l_{\theta}(\varphi_{l,a}) - E_1\parallel\varphi_{l,a}\parallel^2
=\parallel\nabla'\varphi_{l,a}\parallel^2+\parallel\dot{\theta}\partial_{\tau}\varphi_{l,a}+\partial_s
\varphi_{l,a}\parallel^2-E_1\parallel\varphi_{l,a}\parallel^2.
\end{equation}
Evidently the first term on the r.h.s of (\ref{eq353}) is zero.  For
the second term on the r.h.s of (\ref{eq353}) we get,
\begin{eqnarray}\label{eq393}
\nonumber\parallel\dot{\theta}\partial_{\tau}\varphi_{l,a}+\partial_s\varphi_{l,a}\parallel^2
=\parallel \partial_s{ \varphi}_{l,a}\parallel^2 =
\frac{20}{l}\mid\omega\mid.
\end{eqnarray}
 Then
\begin{equation}\label{eq413}
q^l_{\theta}(\varphi_{l,a}) - E_1\parallel\varphi_{l,a}\parallel^2
=\mid\omega\mid(\frac{20}{l}-\frac{13l}{15}E_1)
\end{equation}
and thus if $l \geq l_{min}:=\sqrt{\frac{300}{13E_1}}$ we have
$q^l_{\theta}(\varphi_{l,a}) - E_1\parallel\varphi_{l,a}\parallel^2
\leq0$ $\qed$
\medskip

\subsection{Proof of the Theorem \ref{resultat} $i)$}

Using the same notation as in the Theorem \ref{resultat} $i)$, then
$H^{\mathcal N}_\theta \leq   H^l_\theta $. Moreover these operators
have the same essential spectrum, then by the min-max principle the
assertion follows.
\section{Absence of bound state}

In this section   we want to prove the second part of the Theorem
\ref{resultat}. Denote by $\theta _m = \inf (\rm{supp}(\dot
\theta))$, $\theta _M = \sup($supp$(\dot \theta))$ and
 $L= \theta _M- \theta_m$.  Here $L>0$.  We first consider the case where the Neumann window is an annulus,  $\mathcal{A}_a(l)= I_a(l)\times \omega$.

\begin{proposition}\label{maintheorem}  Suppose
$\mathcal{A}_{a}(l)$ is  such that  $ a \geq \theta _M $. Assume that conditions of the Theorem \ref{resultat} ii) hold. Then
 there exists $d_{max}:=d_{max}(\omega,\theta)>0$, such that for all $0<d \leq d_{max}$ there exists $l_{max}(d,\theta,\omega)>0$ such as for
 all
 $0<l \leq l_{max}$ we have
\begin{equation}
\sigma_{d}(H_\theta^l)=\emptyset.
\end{equation}
\end{proposition}

\remark the case where $l+a \leq \theta _m$ follows  from  same
arguments  developed below.

This   proof is based  on the fact that under conditions of the
Proposition \ref{maintheorem},  for every $\psi \in D(q^l_{\theta})$
it  holds,

\begin{equation} \label{fin}
Q(\psi):= q^l_{\theta}(\psi) - E_1\parallel\psi\parallel^2 \geq 0.
\end{equation}
  The proof of \eqref{fin} involves  several  steps.

\subsection{A local Hardy inequality} The aim of this paragraph is to
show  a   Hardy  type inequality  needed  for the proof of the
Proposition \ref{maintheorem}. It  is the first step of the  proof
of \eqref{fin}. Let  $g$ be the following function
\begin{equation}\label{g}
 g(s):=\left\{
              \begin{array}{ll}
                0, & \hbox{ on} \; \; I_a(l); \\
                E_1, & \hbox{elsewhere.}
              \end{array}
            \right.
\end{equation}
Choose  $p \in (\theta_m, \theta_M)$  s.t. $\dot \theta (p) \not=0$
and let
  \begin{equation}\label{rho}
             \rho(s):=\left\{\begin{array}{ll}
                                    \frac{1}{1+(s-p)^2}, & \hbox{ on $(-\infty,p];$ } \\
                                    0, & \hbox{ elsewhere.}
                                  \end{array}
                                \right.
\end{equation}

\begin{proposition}\label{hardytheorem}
Under  same conditions of the Proposition \ref{maintheorem},
then
 there exists   a constant $C>0$ depending on $p$ and $\omega$ and
$\dot{\theta}$ such that for any $\ \psi\in D(q^l_\theta)$,
\begin{equation}\label{hardy1}
    \parallel\nabla'\psi\parallel^2+\parallel\dot{\theta}\partial_{\tau}\psi+\partial_s\psi\parallel^2-\int_{\Omega_{0}}g(s)\mid\psi\mid^2dsdt\geq C \int_{\Omega_0}\rho(s)\mid\psi\mid^2dsdt.
\end{equation}
\end{proposition}

 We first show  the following lemma. Denote by  $\Omega_{p}:=(-\infty,p)\times\omega$.
\begin{lemma}\label{Lemmap}
Under same conditions of the Proposition \ref{hardytheorem}. Then
for any $\psi\in D(q^l_{{\theta}})$ we have
\begin{equation}\label{un}
\int_{\Omega_p}\mid\nabla'\psi\mid^2
+\mid\dot{\theta}\partial_{\tau}\psi+\partial_s\psi\mid^2-E_1\mid\psi\mid^2dsdt\geq
C\int_{\Omega_{p}} \rho(s) \mid\psi\mid^2 dsdt.
\end{equation}

\end{lemma}

In the following we will  use notations suggested in \cite{EKK}.
 For $ A \subset \mathbb{R} $ denote by  $\chi_A$  the characteristic function of $ A \times \omega$. Let $\psi\in D(q_{\theta}^l)$ and define,
\begin{align}\label{Inotations}
\nonumber q_1^A(\psi)&:=\parallel \chi_A
\nabla'\psi\parallel^2-E_1\parallel \chi_A \psi\parallel^2,\quad
\ q_2^A(\psi):=\parallel  \chi_A \partial_s\psi\parallel^2,\\
q_3^A(\psi)&:=\parallel \chi_A
\dot{\theta}\partial_\tau\psi\parallel^2,\quad
\quad\quad\quad\quad\quad\quad q_{2,3}^A(\psi):=2 \Re(\partial_s\psi,
\chi_A \dot{\theta}\partial_\tau\psi).
\end{align}
Denote also  by  $Q^A (\psi)= q_1^A(\psi)+ q_2^A(\psi)+ q_3^A(\psi)
+ q_{2,3}^A(\psi)$. Here and hereafter  we often  use the fact that
for any $\psi \in D(q_{\theta}^l)$
 \begin{equation}\label{utile}
 q_1^A (\psi)   \geq 0,
\end{equation}
for  every $ A \subset \mathbb R$ such that  $ A \cap I_a(l) =
\emptyset$. \proof Choose $r>0$ such that $\dot{\theta}(s)\neq0$
for any $ s\in[p-r,p]$. Let $f$  be the
 following function:
\begin{equation}f(s):=\left\{
                \begin{array}{ll}
                  0, & \hbox{on $(p,\infty)$;} \\
                  \frac{p-s}{r}, & \hbox{on $(p-r,p]$;} \\
                  1, & \hbox{elsewhere.}
                \end{array}
              \right.
\end{equation}
 For  any $\psi\in D(q^l_{\theta})$,  simple estimates lead to:
\begin{eqnarray}\label{lem2proofex1}
  \int_{\Omega_p}\frac{\mid\psi(s,t)\mid^2}{1+(s-p)^2}dsdt
  &=&\int_{\Omega_{p}}\frac{\mid\psi(s,t)f(s)+(1-f(s))\psi(s,t)\mid^2}{1+(s-p)^2}dsdt  \\
 &\leq&\nonumber
   2\Big(\int_{\Omega_{p}}\frac{\mid
     f(s)\psi(s,t)\mid^2}{(s-p)^2}dsdt+\Vert \chi_{(p-r,p)} \psi \Vert^2\Big).
\end{eqnarray}

 Since     $f(p)\psi(p,.)=0$, we can use the usual Hardy inequality   (see e.g. \cite{Ha}), then  we get,
\begin{equation} \label{expressiontoestimate1}
\int_{\Omega_p}\frac{\mid\psi(s,t)\mid^2}{1+(s-p)^2}dsdt \leq 8
q_2^{(-\infty,p)}(f\psi)+2\Vert \chi_{(p-r,p)} \psi \Vert^2.
  \end{equation}
Note that with our choice $ [p-r,p] \cap [a, a+l] = \emptyset $.
Hence  to estimate the  second  term on the r.h.s of
\eqref{expressiontoestimate1} we    use the  Theorem $6.5$  of
\cite{Kre}, then there exists $ \lambda_0 =
\lambda_0(\dot{\theta},p,r ) >0$ s.t.  for any $\psi\in
D(q_{\theta}^l)$ we have
\begin{equation}\label{secondtermesptimation}
  \Vert \chi_{(p-r,p)} \psi \Vert^2 \leq
\frac{1}{\lambda_0} Q^{(p-r,p)}(\psi)  \leq \frac{1}{\lambda_0}
Q^{(-\infty,p)}(\psi).
\end{equation}
We now  want to estimate the first  term on the right hand side of
(\ref{expressiontoestimate1}).  We have

\begin{equation}\label{eq16t}
q_2^{(-\infty,p)}(f\psi)= \int_{\Omega_p}\mid\partial_s (f
\psi)\mid^2dsdt= q_2^{(-\infty,\theta_m)}(f\psi)+
q_2^{(\theta_m,p)}(f\psi).
\end{equation}
 Evidently  since $ \dot \theta =0 $ and $f=1$ in $(-\infty, \theta_m)$, from \eqref{utile}, we have
\begin{equation}\label{q21}
q_2^{(-\infty,\theta_m)}(f\psi) \leq  Q^{(-\infty,\theta_m)} (\psi).
\end{equation}
In the other hand since   $f(p)\psi(p,.) =0$, we can apply the Lemma
\ref{app}  of the Appendix. So for any $0< \alpha < 1$  there exists
$\gamma_{\alpha,1}>0$ such that
\begin{equation}\label{lem2app}
\mid q^{ (\theta_m,p)}_{2,3} (f\psi)  \mid  \leq
\gamma_{\alpha,1}q_1^{(\theta_m,p)}(f\psi)+  \alpha
q^{(\theta_m,p)}_{2}(f\psi)+q^{(\theta_m,p)}_3(f\psi).
\end{equation}

Let $\gamma:=\max(1,\gamma_{\alpha,1})$.  Then
\begin{eqnarray}\label{lem2app2}
 \gamma^{-1}\mid q^{(\theta_m,p)}_{2,3}(f\psi)\mid\leq q_1^{(\theta_m,p)}(f\psi) + \alpha
\gamma^{-1}q^{(\theta_m,p)}_{2}(f\psi)+\gamma^{-1}q^{(\theta_m,p)}_3(f\psi).
\end{eqnarray}
Hence  with the decomposition,
$q_{2,3}^{(\theta_m,p)}=\gamma^{-1}q_{2,3}^{(\theta_m,p)}+(1-\gamma^{-1})q_{2,3}^{(\theta_m,p)}$
and
 (\ref{lem2app2}) we have,
\begin{align}\label{estimatewithI}
  Q^{(\theta_m,p)} (f\psi)&\geq
(1-\gamma^{-1})\Big(q_{2}^{(\theta_m,p)}(f\psi) +
q^{(\theta_m,p)}_{2,3}(f\psi)+q_{3}^{(\theta_m,p)}(f\psi)\Big) \\ \nonumber
& +\gamma^{-1}(1-\alpha)q^{(\theta_m,p)}_{2}(f\psi)
\end{align}
 and since  $q^{(\theta_m,p)}_{3}+q^{(\theta_m,p)}_{2,3}+q_{2}^{(\theta_m,p)}\geq0$,  we
arrive at,
\begin{equation}\label{firstpart}
 q_{2}^{(\theta_m,p)}(f\psi) \leq\frac{\gamma}{(1-\alpha)}Q ^{(\theta_m,p)}(f\psi).
\end{equation}
Now by using that,   $ q^{(\theta_m,p)}_1(f\psi) \leq
q^{(\theta_m,p)}_1(\psi) $,
$$ \Vert \chi_{(\theta_m,p)}(\partial_s+\dot{\theta}\partial_{\tau})(f\psi)\Vert^2 \leq  2( \Vert\chi_{(\theta_m,p)}(\partial_s+\dot{\theta}\partial_{\tau})\psi\Vert^2 + \frac{1}{r^2} \Vert \chi_{(p-r,p)} \psi \Vert^2)$$
and  \eqref{secondtermesptimation}, we get,
\begin{equation} \label{bri1}
q_{2}^{(\theta_m,p)}(f\psi) \leq
\frac{2\gamma}{(1-\alpha)}(Q^{(\theta_m,p)}(\psi)  +
\frac{1}{\lambda_0r^2}Q^{(p-r,p)}(\psi) ) \leq  c'
Q^{(\theta_m,p)}(\psi)
\end{equation}
 with $c'=   \frac{2\gamma}{(1-\alpha)} (1 +\frac{1}{\lambda_0r^2})$. Then \eqref{q21} and \eqref{bri1} imply
\begin{equation} \label{secondpart}
q_{2}^{(-\infty,p)}(f\psi) \leq (1+c')  Q^{(-\infty,p)}(\psi).
\end{equation}
Hence \eqref {secondpart} and  \eqref {secondtermesptimation}  prove
the lemma  with
\begin{equation}\label{C}
 C^{-1}= 8\big(1+c' \big) +  \frac{2}{\lambda_0}.
\end{equation}
\qed
\bigskip

{ \it Proof of the proposition \ref{hardytheorem}.}  To prove the
proposition we note that \label{lemalpha1alpha2}
 for any $\psi\in D(q_\theta^{l})$ and  for $ p'
\in\mathbb{R}$ we have
\begin{equation}\label{lem1ex}
   \int_{\omega}
   \int_{p'}^{\infty}\mid\nabla'\psi\mid^2+\mid\dot{\theta}\partial_{\tau}\psi+\partial_s\psi\mid^2dsdt\geq\int_{\omega}\int_{p'}^{\infty}g(s)\mid\psi\mid^2dsdt.
\end{equation}

Then    \eqref{lem1ex}  with  $p'=p$ and Lemma \ref{Lemmap} imply
\eqref{hardy1}. \qed

\subsection{Reduction to a one dimensional problem}

We  now want  to prove the following result,
\begin{proposition}\label{proppositivity}  Under  conditions of the Proposition \ref{maintheorem},    then a sufficient condition in order  to get  \eqref{fin} is  given by
\begin{equation}\label{consuf}
\int_{\mathbb{R}}\mid \psi'(s)\mid^2+2C\rho(s)\mid \psi(s)\mid^2ds-
4E_1\int_{a}^{a+l}\mid \psi(s)\mid^2ds \geq 0,  \;\; \;
\end{equation}
for any $ \psi\in \mathcal{H}^1(\mathbb{R})$ where the constant $C$
is defined  in \eqref{C}.
\end{proposition}

 \remark This proposition means that  the positivity  needed here is given by  the positivity of  the effective   one
dimensional Schr\"{o}dinger operator  on $\mathrm{L}^2(\mathbb{R})$,
\begin{equation}\label{schrodingerop}
-\frac{d^2}{ds^2}+2C\rho(s)-4E_1{ \bf 1}_{I_a(l)}.
\end{equation}
where ${ \bf 1}_{I_a(l)}$ is the characteristic function of
$I_a(l)$.


\proof  Evidently we have
\begin{equation}\label{pssi}
   Q(\psi) =
   \frac{1}{2}\big(Q(\psi)-\int_{\Omega_0}(E_1-g(s))\mid\psi\mid^2dsdt +q_\theta^l(\psi)-\int_{\Omega_0}g(s)\mid\psi\mid^2dsdt \big),
\end{equation}
where $g$ is defined in (\ref{g}).  By using  (\ref{hardy1}), then
\begin{equation}\label{reduction1}
Q(\psi)\geq\frac{1}{2}\Big(q_\theta^l(\psi)-E_1\parallel\psi\parallel^2 + \;  C\int_{\Omega_0}\rho(s)\mid\psi\mid^2dsdt-E_1\parallel\chi_{(a,a+l)}\psi\parallel^2\Big)
\end{equation}
 Rewrite the expression of $q_\theta^l$ given
by (\ref{qatheta}) as follows:
\begin{equation}\label{expanded}
q_\theta^l(\psi)=\parallel\nabla'\psi\parallel^2+\parallel\partial_s\psi\parallel^2+\parallel\dot{\theta}\partial_\tau\psi\parallel^2+2 \Re(\partial_s\psi,\dot{\theta}\partial_\tau\psi).
\end{equation}
 We   estimate  the  last  term of the r.h.s. of \eqref{expanded}.
  By using  the formula \eqref {q231estimate}  of the Appendix,
 \begin{equation}\label{young}
\mid q_{2,3}(\psi)\mid= \mid q_{2,3}^{(\theta_m,\theta_M)}(\psi)
\mid \leq\gamma_{\frac{1}{2},\frac{1}{2}} q_1^{(\theta_m,\theta_M)}
(\psi) +\frac{1}{2}q_2^{(\theta_m,\theta_M)}(\psi)
+\frac{1}{2}q_3^{(\theta_m,\theta_M)}(\psi) \end{equation} where
\begin{equation}\label{gamma}
\gamma_{\frac{1}{2},\frac{1}{2}}:=\widetilde{\gamma}_{\frac{1}{2},\frac{1}{2}}+4d^2\parallel\dot{\theta}\parallel_{\infty}^2
 \end{equation} with
$\widetilde{\gamma}_{\frac{1}{2},\frac{1}{2}}:=max\Big\{\frac{d\parallel\dot{\theta}\parallel_{\infty}\parallel\ddot{\theta}\parallel_{\infty}\sqrt{f(L)}}{\dot{\theta_0}\sqrt{\lambda}},\frac{d^2\parallel\ddot{\theta}\parallel^2_{\infty}f(L)}{\lambda\dot{\theta_0}^2},
2d^2\parallel\ddot{\theta}\parallel^2_{\infty}f(L)\Big\}$ for some
constant $\lambda >0$ depending only on the section $ \omega$ and
$f(L):=\max\{2+\frac{16L^2}{r^2},4L^2\}$.
\medskip

 Hence  (\ref{expanded}) together with  (\ref{young})
 give:
\begin{equation}\label{expanded22}
q_\theta^l(\psi) \geq
\parallel\nabla'\psi\parallel^2+\frac{1}{2}\parallel\partial_s\psi\parallel^2+\frac{1}{2}\parallel\dot{\theta}\partial_\tau\psi\parallel^2-\gamma_{\frac{1}{2},\frac{1}{2}}q_1^{(\theta_m,\theta_M)}(\psi).
\end{equation}
In view of  \eqref{utile}   we have

$$\parallel \nabla'\psi\parallel^2-E_1\parallel\psi \Vert^2 \geq  q_1^ {(\theta_m,\theta_M)}(\psi) + q_1^{I_a(l)}(\psi) \geq q_1^ {(\theta_m,\theta_M)}(\psi) -E_1 \Vert \chi_{(a, a+l)} \psi \Vert^2 .$$
 Thus this last inequality together with  (\ref{expanded22}) in (\ref{reduction1})  give

\begin{eqnarray}\label{reduction22}
 \nonumber Q(\psi)
  &\geq&\frac{1}{2}\Big(\frac{1}{2}\parallel\partial_s\psi\parallel^2+\frac{1}{2}\parallel\dot{\theta}\partial_{\tau}\psi\parallel^2+ C \int_{\Omega_0}\rho(s)\mid\psi\mid^2dsdt-2E_1\parallel\chi_{(a,l+a)}\psi\parallel^2\\
\nonumber
&+&(1-\gamma_{\frac{1}{2},\frac{1}{2}})q_1^{(\theta_m\theta_M)}(\psi)\Big).
\end{eqnarray}
Now if  $0<d \leq d_{max}$ then  $\gamma_{\frac{1}{2},\frac{1}{2 }} \leq 1$
so  the Proposition \ref {proppositivity}  follows.

\subsection{The one dimensional Schr\"{o}dinger operator}
 In this part, under our  conditions,
we want to  show   that the   one dimensional Schr\"{o}dinger
operator \eqref{schrodingerop}  is a positive operator.  In view of
the Proposition \ref{proppositivity} this will    imply the
Proposition \ref{maintheorem}. Here we follow a similar  strategy
as  in \cite{BEK}.
\begin{proposition}\label{onedProb}
for all $\varphi\in \mathcal{H}^1(\mathbb{R})$, then there exists
$l_{max}>0$ such  that  for any $ \ 0< l \leq l_{max}$ we have
\begin{equation}\label{schrodinger1d}
\int_{\mathbb{R}}\mid\varphi'(s)\mid^2+2C\rho(s)\mid\varphi(s)\mid^2ds\geq4E_1\int_{I_a(l)}\mid\varphi(s)\mid^2ds.
\end{equation}
\end{proposition}

 \proof
 Introduce the
following function:
\begin{equation}\Phi(s):=\left\{
  \begin{array}{ll}
    (\frac{\pi}{2}+\arctan{(s-p)}), & \hbox{ $if \ s<p$;} \\
    \frac{\pi }{2}, & \hbox{$if \ s\geq p$.}
  \end{array}
\right.
\end{equation}
 where $p$ is the same real  number  as  in \eqref{rho}. So clearly $\Phi'=\rho$. For any $t\in I_a(l)$ and
 $\varphi \in \mathcal{H}^1(\mathbb{R})$, we  have:
\begin{eqnarray}\label{phifirst}
  \frac{\pi }{2}\varphi(t) = \Phi(t)\varphi(t)
   \nonumber&=& \int_{-\infty}^t(\Phi(s)\varphi(s))'ds \\
   &=&\int_{-\infty}^t \rho(s)\varphi(s)ds+\int_{-\infty}^t\Phi(s)\varphi'(s)ds
\end{eqnarray}
and since $\rho(s)=0$ for any $s\in(p,\infty)$, we get,
\begin{equation}\label{phisecond}
\frac{\pi }{2}
\varphi(t)=\int_{-\infty}^{p}\rho(s)\varphi(s)ds+\int_{-\infty}^{t}\Phi(s)\varphi'(s)ds.
\end{equation}
Then some straightforward estimates lead to,
\begin{eqnarray}\label{phithird}
\frac{\pi^2}{4}\varphi^2(t)
&\leq&2\Big((\int_{-\infty}^{p}\rho(s)\varphi(s)ds)^2+(\int_{-\infty}^t\Phi(s)\varphi'(s)ds)^2\Big)\\
&\leq
&2\Big(\int_{-\infty}^{p}\rho(s)ds\int_{-\infty}^{p}\rho(s)\varphi^2(s)ds+\int_{-\infty}^t\Phi^2(s)ds\int_{-\infty}^t\varphi'^2(s)ds\Big)
\nonumber.
\end{eqnarray}

By  direct calculation $ \int_{-\infty}^{p}\rho(s)ds =\frac{\pi }{2}
$ and $\int_{-\infty}^{p}\Phi^2(s)ds+\int_{p}^{t}\Phi^2(s)ds=
\pi\ln{2}+\frac{\pi^2}{4}(t-p)$. Hence we get,
\begin{equation}\label{phi5}
\mid \varphi(t)\mid^2\leq\frac{4}{\pi }\int_{\mathbb{R}}\rho(s)
\varphi^2(s)ds+\Big(\frac{8\ln{2}}{\pi}+2(t-p)\Big)\int_{\mathbb{R}}\mid
\varphi'(s)\mid^2ds
\end{equation}
We integrate  both sides of (\ref{phi5}) over $I_a(l)$, then
\begin{eqnarray}
\nonumber \int_{I_a(l)}\mid \varphi(t)\mid^2 dt&\leq&\frac{4l}{\pi
}\int_{\mathbb{R}}\rho(s)
\varphi^2(s)ds+\Big((\frac{8\ln{2}}{\pi}+2(a-p))l+l^2\Big)\int_{\mathbb{R}}\mid
\varphi'(s)\mid^2ds\\
&\leq& c''\int_{\mathbb{R}}2C\rho(s)
\varphi^2(s)+ \mid \varphi'(s)\mid^2ds
\nonumber
\end{eqnarray}
 where $c'' = 2l(\frac{1}{\pi C}+\frac{4\ln{2}}{\pi}+a-p)+l^2$. Finally we get,
\begin{equation}
4E_1\int_{a}^{l+a}\mid \varphi(t)\mid^2dt\leq
4E_1c'' \int_{\mathbb{R}}2C\rho(s)\mid \varphi(s)\mid^2+\mid
\varphi'(s)\mid^2ds.
\end{equation}
So  choose $0< l \leq  l_{max}$  with
$$ l_{max}:= -(\frac{1}{\pi C}+\frac{4\ln{2}}{\pi}+a-p)+\sqrt{(\frac{1}{\pi C}+\frac{4\ln{2}}{\pi}+a-p)^2 + (4E_1)^{-1}}$$
   then  $4E_1c'' \leq 1$ and  the proposition
\ref{onedProb} follows. \qed

\subsection{ proof of the Theorem \ref{resultat} $ii)$}
  Under   assumptions of the  Theorem \ref{resultat} $ii)$  $H^{\mathcal N}_\theta \geq   H^l_\theta $. These two operators  have the same essential spectrum so the  Theorem \ref{resultat} $ii)$  is proved   by applying the Proposition \ref{maintheorem}  and the min-max principle.

\section{Appendix}

In this appendix we give a slight extension of the lemma $3$ of
\cite{EKK}   which  states that  under our conditions, for all
$\psi\in D(q_\theta^l)$   we have  for any $\alpha,\beta >0$ there
exists $\gamma_{\alpha,\beta}>0$ such that:
\begin{equation}\label{q231estimate}
\mid q_{2,3}(\psi)\mid\leq \gamma_{\alpha,\beta}q_1(\psi)+\alpha
q_2(\psi)+\beta q_3(\psi).
\end{equation}
Then we have
\begin{lemma} \label{app} Let  $p \in (\theta_m, \theta_M)$.
For all $\psi\in D(q_\theta^l)$  such that $\psi(p,.)=0$, then
for any $ \alpha,\beta >0$ there exists $\gamma_{\alpha,\beta}>0$
such that:
\begin{equation}\label{q23estimate}
\mid q_{2,3}^{(\theta_m,p)}(\psi)\mid\leq
\gamma_{\alpha,\beta}q_1^{(\theta_m,p)}(\psi)+\alpha
q_2^{(\theta_m,p)}(\psi)+\beta q_3^{(\theta_m,p)}(\psi).
\end{equation}
\end{lemma}
\proof Let $\psi\in D(q_\theta^l)$ such that   $\psi(p,.)=0$.  Then
$\psi\in \mathcal{H}_0^1(\Omega_p)$.
 We know that   we may  first   consider  vectors  $\psi(s,t)=\chi_1(t)\phi(s,t)$, where $\phi\in C_0^\infty
(\Omega_p)$.  For such a vector $\psi$ we have,
\begin{eqnarray}  q_1^{(\theta_m,p)}(\psi)&=&\parallel\chi_{(\theta_m,p)}\chi_1\nabla'\phi\parallel^2,\quad
q_2^{(\theta_m,p)}(\psi)=\parallel\chi_{(\theta_m,p)}\chi_1\partial_s\phi\parallel^2\\ \nonumber
q_3^{(\theta_m,p)}(\psi)&=&\parallel\chi_{(\theta_m,p)}\dot{\theta}(\chi_1\partial_\tau\phi+\phi\partial_\tau\chi_1)\parallel^2
\end{eqnarray}
and
\begin{equation}\label{q23Lp}
q_{2,3}^{(\theta_m,p)}(\psi)=2(\dot{\theta}\chi_{(\theta_m,p)}
\chi_1\partial_\tau\phi,\chi_1\partial_s\phi)+2(\dot{\theta}\chi_{(\theta_m,p)}\phi\partial_\tau\chi_1,\chi_1\partial_s\phi)
\end{equation}
By using simple estimates  the first term on the r.h.s of
(\ref{q23Lp}) is  estimated as :
\begin{eqnarray} \nonumber
\mid2(\dot{\theta}\chi_{(\theta_m,p)}\chi_1\partial_\tau\phi,\chi_1\partial_s\phi)\mid\leq 2\parallel\dot{\theta}\parallel_\infty\parallel\chi_{(\theta_m,p)}\chi_1\nabla'\phi\parallel\parallel\chi_{(\theta_m,p)}\chi_1\partial_s\phi\parallel
\end{eqnarray}
then
\begin{equation}  \label{n1}
\mid2(\dot{\theta}\chi_{(\theta_m,p)}\chi_1\partial_\tau\phi,\chi_1\partial_s\phi)\mid
 \leq  c_1
q_1^{(\theta_m,p)}(\psi)+\frac{\alpha}{2}q_2^{(\theta_m,p)}(\psi),
\end{equation}
where
$c_1:=\frac{2}{\alpha}d^2\parallel\dot{\theta}\parallel^2_\infty$
and $\alpha>0$.\\
  Integrating by parts twice and  using  the
fact that $\dot{\theta}(\theta_m)=\phi(p,.)=0$,  the second
term of the r.h.s of (\ref{q23Lp})  is written as
\begin{equation}\label{2estimate}
2(\dot{\theta}\chi_{(\theta_m,p)}\phi\partial_\tau\chi_1,\chi_1\partial_s\phi)=(\chi_{(\theta_m,p)}\ddot{\theta}\phi\chi_1,\chi_1\partial_\tau\phi).
\end{equation}
Then the Cauchy Schwartz inequality implies,
\begin{equation}
\mid(\chi_{(\theta_m,p)}\ddot{\theta}\phi\chi_1,\chi_1\partial_\tau\phi)\mid^2\leq
d^2\parallel\ddot{\theta}\parallel^2_\infty q_1^{(\theta_m,p)}
\parallel\chi_{(\theta_m,p)}\chi_1\phi\parallel^2.
\end{equation}
Let $p'  \in \mathbb R$ and $ r'>0$ such that  $(p'-r,p')\subset
(\theta_m,p)$
 and  for $ s\in
(p'-r,p')$,   $\mid\dot{\theta}(s)\mid\geq \dot{\theta}_0$  for some
$\dot{\theta}_0>0$. As in the proof of the Lemma $3$ of  \cite{EKK}
we have,
\begin{equation}\label{Lem1EKK}
 \parallel\chi_{(\theta_m,p)}\chi_1\phi\parallel^2\leq c_2\Big(q_2^{(\theta_m,p)}(\psi)+\dot{\theta}_0^{-2}\parallel\chi_{(p'-r,p')}\dot{\theta}\chi_1\phi\parallel^2\Big)
\end{equation}
where $c_2:=\max\Big\{2+16\frac{(p-\theta_m)^2}{r^2},4(p-\theta_m)^2\Big\}$.\\
Moreover, for any $s\in \mathbb{R}$,
$\dot{\theta}(s)\chi_1\phi(s,.)\in\mathcal{H}_0^1(\Omega_p)$, then
by
 using the  Lemma $1$ of \cite{EKK} there exists $\lambda>0$ depending
 on $\omega$ such that :
 \begin{equation}\label{EKKestimate}
 \parallel\chi_{(p'-r,p')}\dot{\theta}\chi_1\phi\parallel^2\leq\parallel\chi_{(\theta_m,p)}\dot{\theta}\chi_1\phi\parallel^2\leq\lambda^{-1}\Big(q_3^{(\theta_m,p)}(\psi) +\parallel\dot{\theta}\parallel^2_\infty q_1^{(\theta_m,p)}(\psi)\Big).
 \end{equation}
Hence  (\ref{Lem1EKK}), (\ref{EKKestimate})   and  (\ref{2estimate})
give
 \begin{equation} \label{n2}
\mid(\chi_{(\theta_m,p)}\ddot{\theta}\phi\chi_1,\chi_1\partial_\tau\phi)\mid^2\leq\Big(c_3
q_1^{(\theta_m,p)}(\psi)+\frac{\alpha}{2} q_2^{(\theta_m,p)}(\psi)
+\beta q_3^{(\theta_m,p)}(\psi)\Big)^2
 \end{equation}
 where $c_3:=\max\Big\{\frac{d\parallel\ddot{\theta}\parallel\parallel\dot{\theta}\parallel_\infty\sqrt{c_2}}{\dot{\theta}_0\sqrt{\lambda}},\frac{d^2\parallel\ddot{\theta}\parallel^2_\infty c_2}{\alpha},\frac{d^2\parallel\ddot{\theta}\parallel^2_\infty c_2}{2\beta\dot{\theta}_0^2\lambda}\Big\}$. Then \eqref{n1} and \eqref{n2}
imply \eqref{q23estimate} with $\gamma_{\alpha,\beta}:= c_1+c_3$.

Note that  we can choose $\chi_1 >0$  on $\omega$. So that
\eqref{q23estimate}  holds  for every $\psi \in
 C_0^\infty
(\Omega_p)$ and by  a density  argument  this   is even true  for
$\psi\in \mathcal{H}_0^1(\Omega_p)$.


\end{document}